\documentclass{amsart}
\usepackage{amssymb,amsmath,amsthm,amsfonts,amsopn,url,bbm}
\usepackage{enumitem, color}
\usepackage{amscd, hyperref}
\theoremstyle{plain}
\newtheorem{thm}{Theorem}

\theoremstyle{definition}
\newtheorem{defn}[thm]{Definition}
\newtheorem*{defn*}{Definition}
\newtheorem{example}[thm]{Example}
\newtheorem{construction}[thm]{Construction}
\theoremstyle{remark}

\newtheorem*{rmk*}{Remark}
\newcommand{\field}[1]{\mathbbm{#1}}
\newcommand{\N}{\field{N}}

\newcommand{\ideal}[1]{\mathfrak{#1}}
\newcommand{\m}{\ideal{m}}

\newcommand{\p}{\ideal{p}}

\newcommand{\ia}{\ideal{a}}
\newcommand{\ib}{\ideal{b}}
\newcommand{\ie}{\ideal{e}}
\newcommand{\ih}{\ideal{h}}
\newcommand{\ic}{\ideal{c}}

\newcommand{\func}[1]{\mathrm{#1} \,}

\newcommand{\Spec}{\func{Spec}}

\DeclareMathOperator{\Ass}{Ass}

\newcommand{\ra}{\rightarrow}

\newcommand{\rjj}[3]{u_{#1}({#2}, {#3})}

\newcommand{\rjjm}[3]{u^-_{#1}({#2}, {#3})}

\newcommand{\rjjp}[3]{u^+_{#1}({#2}, {#3})}

\DeclareMathOperator{\len}{\lambda}

\DeclareMathOperator{\lt}{\mathtt{lt}}
\DeclareMathOperator{\lm}{\mathtt{lm}}
\DeclareMathOperator{\lcm}{\mathtt{lcm}}

\title{A computation concerning relative Hilbert-Kunz multiplicities}

\author{Neil Epstein}
  \address{Department of Mathematical Sciences \\ George Mason University \\ Fairfax, VA  22030}
  \email{nepstei2@gmu.edu}

\author{Yongwei Yao}
\address{
Department of Math and Statistics\\
Georgia State University\\
30 Pryor St., Atlanta, GA 30303}
\email{yyao@gsu.edu}
\thanks{The second author was partially supported by the National Science Foundation DMS-0700554}


\date{\today}
\begin{document}


\maketitle

In \cite{nmeYao-HK}, we developed several methods designed to provide numerical critera for when a nested pair of submodules of a finitely generated module admit the same tight closure.  For the purposes of this note, it is enough to consider pairs of ideals $J\subseteq I$ such that $\len(I/J)=\infty$..  Hence, we choose to state the definitions and results from that paper in terms of ideal containment.

\begin{defn}\label{def:rjj}
Let $J \subseteq I$ be ideals of a local prime characteristic Noetherian ring $(R,\m)$ of dimension $d$.  Then their \emph{relative multiplicity} is \[
\rjjp R J I := \limsup_{q \ra \infty} \frac{\len(\Gamma_\m(I^{[q]} / J^{[q]}))}{q^d}.
\]
(resp. \[
\rjjm R J I := \liminf_{q \ra \infty} \frac{\len(\Gamma_\m(I^{[q]} / J^{[q]}))}{q^d}).
\]
If these are equal (\emph{i.e.}, the limit is well-defined), then the common number is written $\rjj RJI$.
\end{defn}

Recall the following Theorem, specialized to the case of ideal containment.

\begin{thm}[{\cite[part of Theorem 2.4]{nmeYao-HK}}]\label{thm:zerotc}
Let $R$ be a Noetherian ring, and let $J \subseteq I$ be ideals.  Suppose that $R$ contains a completely stable weak test element $c$, and that $\widehat{R_\p}$ is equidimensional for all $\p \in \Spec R$.  If $\rjjm{R_\p}{J_\p}{I_\p}=0$ for all $\p \in \Spec R$, then $I \subseteq J^*$.
\end{thm}

Seeking a converse to Theorem~\ref{thm:zerotc}, let $J \subseteq I$ be ideals with the same tight closure.  In \cite[Proposition 3.1, Theorem 3.4, and Theorem 3.5]{nmeYao-HK}, we gave several criteria under which a converse to Theorem~\ref{thm:zerotc} holds.  An analysis of the ideas surrounding 
\cite[Proposition 3.1]{nmeYao-HK} yields the following observation: \emph{The critical situation occurs when there exist prime ideals $\p \subsetneq \m$ such that $\p, \m \in \Ass_R(I^{[q]}/J^{[q]})$ for infinitely many values of $q$.}  One may ask whether this can happen.  For instance, in 
\cite[Example 2.2]{nmeYao-HK}, the critical situation does not occur for the ideals $J \subseteq I$ in $R$ unless it already was an issue for the ideals $\ib \subseteq \ia$ in $A$.  Indeed, for each $q$, there is a bijective correspondence between the sets $\Ass_A(\ia^{[q]}/\ib^{[q]})$ and $\Ass_R(I^{[q]}/J^{[q]})$, given by $\p \mapsto \p R$.

However, the situation outlined above can happen, as shown below.  Moreover, the expected converse to Theorem~\ref{thm:zerotc} holds, at least in the given example.  Note that 
the example below does not appear to arise as one of the special cases delineated in \cite[\S 3]{nmeYao-HK}.  Therefore, we had to use computational methods.

Before we get to the specific characteristic $p$ situation, we give a somewhat more general construction, which works over any field, and may be of independent interest.  As we will be using Gr\"obner basis techniques, we set some notations and recall some facts:

\begin{defn}
Let $A$ be a polynomial ring, $>$ a monomial order, and $f\in A \setminus \{0\}$.  The expressions $\lt(f)$ and $\lm(f)$ denote, respectively, the leading term and the leading monomial of $f$ with respect to the given order.

Given two elements $f,g \in A \setminus \{0\}$, the \emph{S-polynomial} of $f$ and $g$ is given by \[
S(f,g) := \frac{\lcm(\lt(f), \lt(g))}{\lt(f)} \cdot f - \frac{\lcm(\lt(f), \lt(g))}{\lt(g)} \cdot g,
\]
where $\lcm$ means the least common multiple.
\end{defn}

The following theorem is a slightly nonstandard (albeit well-established) form of the \emph{Buchberger criterion}:

\begin{thm}\label{thm:Bucrit}\cite[Theorem 2.9.3]{CLO-book1e3}
Let $A$ be a polynomial ring over a field, let $>$ be a monomial order, and let $G = \{g_1, \dotsc, g_n\}$ be a finite subset of $A$.  Then $G$ is a Gr\"obner basis if and only if there exist elements $a_{ijk} \in A$ such that for each pair $(j,k)$ with $1\leq j <k \leq n$, we have \[
S(g_j,g_k) = \sum_{i=1}^n a_{ijk} g_i,
\]
in such a way that for each nonzero $a_{ijk}$, we have $\lm(S(g_j, g_k)) \geq \lm (a_{ijk} g_i)$ with respect to the given monomial order.
\end{thm}

\begin{thm}\label{thm:elim} \cite[Theorem 4.3.11 and the discussion which follows]{CLO-book1e3}
Let $A$ be a polynomial ring over a field $k$, let $I$ be an ideal of $A$ and $0 \neq u\in A$.  Let $r$ be an indeterminate over $A$, and let $B = A[r]$ be a polynomial ring, ordered with lexicographic order in such a way that $r>x$ for all variables $x$ of $A$.  Let $\ia := rIB + (1-r)uB \subseteq B$.  Then $\ia \cap A = I \cap (u)$, and if $F$ is a Gr\"obner basis of $\ia$ in $B$, then $\ia \cap A$ is generated by the set of elements of $F$ whose leading terms are not multiples of $r$.
\end{thm}

\begin{construction}\label{cons}
Let $k$ be an arbitrary field, let $m \in \N$ such that $m\geq 4$, and let $n=2m+1$. We also impose the condition that if $p$ is the characteristic of $k$, then $p \nmid m$, which is automatically satisfied if $p=0$. Let $A := k[s,x,y]$, $\m := (s,x,y) \subseteq A$, $g=xy(x-y)(x+y-sy)$, and $\ie := (x^n, y^n, g) \subseteq A$.  Let $f := \sum_{j=2}^{n-1} (-1)^j x^{n+1-j} y^j$.  Let $\ih := \ie + (f)$.  Let $\ib := (x,y)^{n+2}$.  Then we will show the following: \begin{enumerate}
\item $\ib \subseteq \ie$,
\item $sf \in \ie$ (hence, $\ih \subseteq (\ie :s)$),
\item $xf, yf \in \ie$ (hence, $\m \subseteq (\ie :f)$),
\item $f \notin \ie$ (hence, $(\ie:f) \neq A$, so that $\m = (\ie :f)$),
\item $\ih$ is $s$-saturated  (that is, $(\ih : s) = \ih$),
\item $\ie : \m^\infty = \ie : s^\infty = \ih$, and
\item $H^0_\m(A/\ie) \cong A/\m$.
\end{enumerate}

To see (1), take a typical monomial generator $x^i y^j$ of $\ib$.  That is, $i+j = n+2$.  Since $x^n, y^n \in \ie$, we may assume that $1\leq j \leq n-1$, so that $i\geq 3$.  Note that modulo $g$, we have \[
x^3y \equiv s x^2 y^2 - (s-1) x y^3.
\]
Multiplying this by $x^{i-3} y^{j-1}$, we have $x^i y^j \in (x^{i-1} y^{j+1}, x^{i-2} y^{j+2}, g)$.  Then apply induction to obtain $x^i y^j \in (x^2 y^n, x y^{n+1}, g) \subseteq \ie$.

To see (2), note that modulo $g$, we have \[
s x y^2 (x-y) \equiv xy (x^2 - y^2).
\]
Using this congruence, we have:
\begin{align*}
s(f - xy^n) &= s x y^2(x-y)\left(\sum_{j=0}^{m-1} x^{n-3-2j} y^{2j}\right) \\
&\equiv xy(x^2 -y^2)\left(\sum_{j=0}^{m-1} x^{n-3-2j} y^{2j}\right) \\
&= xy(x^{2m} - y^{2m}) = x^{n} y - x y^{n}.
\end{align*}
Thus, $sf \in (x^n, y^n, g) = \ie$, as required.

To see (3), let $t=s-1$.  Modulo $g$, we have the equivalence \[
tx y^2 (x-y) \equiv x^2 y (x-y).
\]
It follows by induction (on $i$) that for all integers $i\geq 1$, $a\geq 1$, and $b\geq i+1$, we have $t^{i} x^a y^b (x-y) \equiv x^{a+i} y^{b-i} (x-y)$ (modulo $g$).  In particular (letting $a=1$ and $b=n$), for all $1\leq i\leq n-1$, we have \[
t^i x y^n(x - y)\equiv x^{i+1}y^{n-i} (x-y).
\]
modulo $g$. Note also that $-x f+ x^{n} y^2 =yf - x^2 y^n = \sum_{j=1}^{m-1} x^{2j+1} y^{n-2j}(x-y)$.  But by the above (since $2(m-1) = n-1$), this latter sum is congruent (modulo $g$) to $y^n \cdot \left(\sum_{j=1}^{m-1} t^{2j} x (x-y) \right)$.  Thus, $-xf, yf \in (x^n, y^n, g) = \ie$, as required.

In order to demonstrate (4), we require the introduction of Gr\"obner bases into the discussion.  From now on, we will use \emph{lexicographic}\footnote{We emphasize here that we are \emph{not} using degree-lexicographic order.  So for instance, in this ordering, we have $s>x^2$.  Indeed, $s>x^{200}$.} order, with $s>x>y$.  We claim that \[
G := \{g, x^n, x^{n-1} y^3, x^{n-2} y^4, \cdots, x^3 y^{n-1}, y^{n} \}
\] is a Gr\"obner basis of $\ie$ with respect to lex order.  First, since the elements of $G$ consists of the generating set $\{g, x^n, y^n\}$ of $\ie$ along with some elements of $\ib$ (an ideal which by (1) is contained in $\ie$), it follows that $G$ is indeed a generating set for $\ie$.  To show that it is a Gr\"obner basis, we shall find $a_{ijk}$ as in Theorem~\ref{thm:Bucrit}.  But since the S-polynomial of a pair of monomials is always 0, we only need to look at the S-polynomials $S(m,g)$ for monomials $m$ of $G$.  In the following list, we represent each S-polynomial in two ways.  First, we write it in lexicographic order, and then we write it in the form given by Theorem~\ref{thm:Bucrit}:  \begin{itemize}
\item $S(y^{n}, g) = -sxy^{n+1} - x^3y^{n-1} + xy^{n+1} = (-sxy)y^{n} - 1(x^3 y^{n-1}) + (xy)y^{n}$.
\item $S(x^3 y^{n-1}, g) = -s x^2 y^{n} - x^4 y^{n-2} + x^2 y^{n} = (-s x^2) y^{n} - 1(x^4 y^{n-2}) + (x^2) y^{n}$.
\item For any $i$ with $4 \leq i \leq n-2$, we have $x^{i-1} y^{n+3-i}, x^{i+1} y^{n+1-i} \in G$.  And \begin{align*}
S(x^i y^{n+2-i}, g) &= - s x^{i-1} y^{n+3-i} - x^{i+1} y^{n+1-i} + x^{i-1} y^{n+3-i} \\
&= (-s+1) x^{i-1} y^{n+3-i} + (-1)x^{i+1}y^{n+1-i}.
\end{align*}
\item $S(x^{n-1} y^3, g) = -s x^{n-2}y^4 - x^n y^2 + x^{n-2} y^4 = (-s+1) x^{n-2} y^4 - (y^2) x^{n}$.
\item $S(x^{n}, g) = -s x^{n-1} y^3 - x^{n+1}y + x^{n-1} y^3 = (-s+1) x^{n-1} y^3 + (xy) x^{n}$.
\end{itemize}
Thus, $G$ is a Gr\"obner basis of $\ie$.  The leading term $x^{n-1} y^2$ of $f$ is manifestly not divisible by any of the leading terms of $G$, which means that the output of the division algorithm of $f$ by $G$ is $f$.  Thus, $f \notin \ie$, as required.

To demonstrate (5), we will use Theorem~\ref{thm:elim}.  Accordingly, let $B := k[r,s,x,y]$, ordered lexicographically with $r>s>x>y$, and consider the ideal $\ia := r\ih B+ (1-r) sB$ of $B$.  We claim that the entries of the following vector comprise a Gr\"obner basis of $\ia$.  Note that it ends with all the elements of $s G \cup \{sf\}$ except for $s x^{n-1}y^3$.
\[
\left[ \begin{matrix}
rs-s\\
r x^{n}\\
c := r x^3y - r x y^3 - s x^2 y^2 + s x y^3 \\
d := m r x^2 y^{n-1} + \sum_{j=1}^{n-3} (-1)^{j-1} j s x^{n-1-j} y^{j+2} \\
r y^{n} \\
-sg = s^2 x^2 y^2 - s^2 xy^3 - s x^3y + s x y^3\\
s x^{n} \\
s f = \sum_{j=2}^{n-1} (-1)^j s x^{n+1-j}y^j \\
s x^{n-2} y^4 \\
s x^{n-3} y^5 \\
\vdots \\
s x^3 y^{n-1} \\
s y^{n}
\end{matrix} \right]
\]
(This is a vector of length $n+5$, and we label the elements $F_0$ through $F_{n+4}$.) First we have to show that the ideal generated by the entries of $F$ is exactly $\ia$.  To see that $\ia \subseteq (F)$, \begin{itemize}
\item $rg = (-x^2 y^2 + x y^3) (rs-s) + 1\cdot c$, and
\item $rf = (x-y) \left(\displaystyle \sum_{j=1}^{m-1} j x^{n-3-2j} y^{2j-1}\right)c + d + mx (-r y^{n} + s y^{n})$.
\end{itemize}
To see that $F \subseteq \ia$, \begin{itemize}
\item $c = 1\cdot (rg) + (-x^2 y^2 + x y^3) (-rs+s)$,
\item $d =(m-1)(-sx+x)(ry^{n}) + (-x+y)\left(\displaystyle \sum_{j=1}^{m-1} j x^{n-3-2j} y^{2j-1}\right)(rg) + 1\cdot (rf)+ \left(\displaystyle \sum_{j=1}^{n-3} (-1)^{j-1} jx^{n-j-1}y^{j+2}\right)(-rs+s)$,
\end{itemize}
and for each element $u \in G \cup \{f\}$, we have $ru \in \ia$, so that \begin{itemize}
\item $su = s\cdot (ru) + u \cdot(-rs+s)$.
\end{itemize}
Thus, $\ia = (F)$.

Taking all the S-polynomials $S_{jk} = S(F_j, F_k)$ such that $j<k$ and $F_j$, $F_k$ are not both monomials (and note that the only non-monomials are $F_i$ for $i=0,2,3,5,7$), we may obtain the following list.  For these choices of $a_{ijk}$, the diligent reader may easily verify the conditions of Theorem~\ref{thm:Bucrit}:
\begin{itemize}
\item $S_{01} = -F_6$
\item $S_{02} = x y^3 F_0 + F_5$
\item $S_{03} = ((-x+y)\sum_{j=1}^{m-1} j x^{n-3-2j} y^{2j-1})F_5  - F_7 + (m-1)(sx -x) F_{n+4}$
\item $S_{04} = -F_{n+4}$
\item $S_{05} = (s x y^3 + x^3y-xy^3)F_0 - F_5$
\item $S_{06} = -F_6$
\item $S_{07} = (\sum_{j=3}^{n-1} (-1)^{j-1} x^{n-j+1}y^j) F_0 - F_7$
\item $S_{0i} = -F_i$, for $8\leq i \leq n+4$
\item $S_{12} = (\sum_{j=1}^{m-1} x^{n-3-2j} y^{2j}) F_2 + xF_4 + F_7 - x F_{n+4}$
\item $S_{13} = (\sum_{j=1}^{n-3} (-1)^j j x^{n-3-j} y^{j+2})F_6$
\item $S_{15} =  (\sum_{j=1}^{n-2} r x^{n-2-j} y^j)F_5  + (rxy+ry^2)F_6 + (rsxy-rx^2 - rxy)F_{n+4}$
\item $S_{17} = - sx^2 F_4+ryF_7$
\item $S_{23} = -mxyF_4 - y F_7 + (\sum_{j=1}^{n-4} (-1)^{j-1} jF_{j+7}) + ((1-m)x^2 + mxy)F_{n+4}$
\item $S_{24} = -x y^2 F_4 + (-x^2 y + x y^2) F_{n+4}$
\item $S_{25} = (s x^2 y^3 - s x y^4 + x^4 y - x^2 y^3) F_0 + (-sy-x)F_5$
\item $S_{26} = -s(\sum_{j=1}^{m-1} x^{n-3-2j}y^{2j})F_2 - sxF_4 - sF_7 + sxF_{n+4}$
\item $S_{27} = (s x^{n-5} y^2 + (-x+y)\sum_{j=2}^{m-1} jsx^{n-3-2j}y^{2j-1})F_2 - s F_3 + (m-1)(rx-sx)F_{n+4}$
\item $S_{2i} =  -sF_{i+1} + (-r+s) F_{i+2}$, for $8\leq i \leq n+1$
\item $S_{2,n+2} = -sF_{n+3} + (-r+s)x^2 F_{n+4}$
\item $S_{2,n+3} = (-rxy-sx^2+sxy)F_{n+4}$
\item $S_{2,n+4} = (-rxy^2 - sx^2y + s xy^2)F_{n+4}$
\item $S_{34} = (\sum_{j=1}^{n-4}(-1)^{j-1} j F_{j+7}) - (n-3) x^2 F_{n+4}$
\item $S_{35} = m(sxy^n + x^3 y^{n-2} - xy^n)F_0 + [s(x-y)(\sum_{j=1}^{m-1} j x^{n-3-2j} y^{2j-1}) - m y^8] F_5 +sF_7 + (m-1) (-s^2 x + s x) F_{n+4}$
\item $S_{36} = (\sum_{j=1}^{n-3} (-1)^{j-1} j x^{n-j-3}y^{j+2}) F_6$
\item $S_{37} = ms(\sum_{j=0}^{n-4} (-1)^j x^{n-2-j}y^j) F_4 + (\sum_{j=1}^{n-4} (-1)^{j-1} jsx^{n-4-j}y^{j+2})F_6 - (n-3)sxy^{n-5} F_8$
\item $S_{3i} = x^{n+3-i} (sy F_7 + s(\sum_{j=1}^{n-4} (-1)^j jF_{j+7}) - sx^2 F_{n+4})$ for $8\leq i\leq n+3$
\item $S_{3,n+4} = (\sum_{j=1}^{n-4} (-1)^{j-1}js F_{j+7}) - (n-3) sx^2 F_{n+4}$
\item $S_{45} = (s^2 xy-sxy) F_4 + r F_{n+3}$
\item $S_{47} = (\sum_{j=1}^{n-3} (-1)^{j-1} s x^{n-1-j} y^j) F_4$
\item $S_{56} = -(\sum_{j=1}^{n-2} x^{n-2-j} y^j)F_5 -(xy+y^2)F_6 + (-sxy+x^2+xy)F_{n+4}$
\item $S_{57} = -(\sum_{j=1}^{m-1} x^{n-3-2j} y^{2j}) F_5 -yF_6+(-sx+x)F_{n+4}$
\item $S_{58} = -yF_7 - F_8 - (s+2)F_9 + (\sum_{j=10}^{n+3} (-1)^{j-1}F_j)+ x^2 F_{n+4}$
\item $S_{5i} = - F_{i-1} + (-s+1)F_{i+1}$ for $9\leq i\leq n+2$
\item $S_{5,n+3} =  -F_{n+2}+ (-s+1)x^2 F_{n+4}$
\item $S_{5,n+4} = - F_{n+3} + (-s+1)xyF_{n+4}$
\item $S_{67} = y F_7 - x^2 F_{n+4}$
\item $S_{7i} = xy^{i-8} ((\sum_{j=9}^{n+3} (-1)^j F_j) + (-x^2+xy)F_{n+4})$ for $8\leq i\leq n+3$
\item $S_{7,n+4} = (\sum_{j=1}^{n-3} (-1)^j x^{n-1-j} y^j)F_{n+4}$
\end{itemize}
Hence, the entries of $F$ give a Gr\"obner basis of $\ia$. By Theorem~\ref{thm:elim}, it follows that the elements of $F$ whose leading term does not involve $r$ forms a generating set for the ideal $\ih \cap (s)$ of $A$.  That is, $\ih \cap (s) = (s y^{n}, s x^3 y^4 (x,y)^{n-5}, s f, s x^{n}, sg)$.  Dividing by $s$, we get $(\ih : s) = (y^{n}, x^3 y^4(x,y)^{n-5}, f, x^{n}, g) = \ih$ (since $x^3 y^4 (x,y)^{n-5} \subseteq (x,y)^{n+2} = \ib \subseteq \ih$), as required. 

To see (6), first note that $\ie : \m^\infty = \ie : s^\infty$, since $\ie$ contains powers of both $x$ and $y$.  But $\ie \subseteq \ih$, so  from (2) and (5), we have $\ih \subseteq (\ie:s) \subseteq (\ih:s) = \ih$, whence all are equalities.  Thus, $(\ie : s^\infty) = (\ih : s^\infty) = \ih$, as required.

Finally, to see (7), it follows from (6) and (4) that \[
H^0_\m(A/\ie) = \frac{\ie : \m^\infty}{\ie} = \frac{\ih}{\ie} = \frac{\ie + (f)}{\ie} \cong \frac{A}{(\ie :f)} = A/\m.
\]
\end{construction}

\begin{example}
Let $p$ be an odd prime number.  Let $k$ be a field of characteristic $p$, and $R := k[s,x,y] / (xy (x-y)(x+y-sy))$.  This is the ring used by Katzman in \cite{KaUAss}, with variable change given by $s=t+1$.

Consider the ideals $J := (x^p, y^p)$ and $I=(x,y)^p$ of $R$.  As shown in Katzman's paper, $J^* = I$.  Now fix a power $q=p^e$ of $p$, $e\geq 1$, and let $n=pq$ in Construction~\ref{cons}.  Let $\ib$, $\ie$, $A$, $\m$, $g$, $\ih$, and $f$ be as in that construction.  The conditions of the construction are satisfied, since $p\geq 3$, whence $n= pq \geq 9$, and $p$ can never divide $(pq-1)/2$.  Then $R=A/(g)$ and $J^{[q]} = \ie/(g) \subseteq R$.  In particular, letting $z$ be the image of $f$ in $R$, we have $z\notin J^{[q]}$.  

However, we claim that $z\in I^{[q]}$.  To see this, it is enough to show (in the ring $R$ -- that is, modulo $g$) that for all $j=2, 3, \dotsc, pq-1$, we have $x^j y^{pq+1-j} \in (x^q, y^q)^p$.   For $j=q, q+1$ this is clear, and for $j\geq q+2$, the assertion follows from the equation $x^j y^{pq+1-j} = (t+1) x^{j-1} y^{pq-j+2} - t x^{j-2} y^{pq-j+3}$, along with induction, showing that in these cases, $x^j y^{pq+1-j} \in (x^q y^{pq-q})$. For $2\leq j \leq q-1$, we have \begin{align*}
x^j y^{pq+1-j} 
&= x^{q+1} y^{pq-q} - \sum_{i=1}^{j-2} x^{q+i} y^{pq-q-i}(y-x) - x^j y^{pq-q+2-j}(x^{q-1}-y^{q-1}) \\
&= x^{q+1} y^{pq-q} - \sum_{i=1}^{j-2} t^i x^{q} y^{pq-q}(y-x) - t^{j-1}x y^{pq-q+1}(x^{q-1}-y^{q-1}) \\
&= x^{q+1} y^{pq-q} - \sum_{i=1}^{j-2} t^i x^{q} y^{pq-q}(y-x) - t^{j-1}x^q y^{pq-q+1}
   + t^{j-1}x y^{pq} \\
& \in (x^q, y^q)^p = I^{[q]}.
\end{align*}

Let $\p := (x,y)$.  We claim that  $\p \in \Ass_R(I^{[q]}/ J^{[q]})$.  Since $\p$ is minimal over $J^{[q]}$, it suffices to show that $I^{[q]}_\p/ J^{[q]}_\p = (I^{[q]} / J^{[q]})_\p \neq 0$.  To do this, it suffices to show that $((x^{pq},y^{pq}) +(g))A_P$ is properly contained in $((x^q, y^q)^p + (g))A_P$, where $P := (x,y) \subseteq A$.  But in the ring $C := L[x,y]$ (where $L:=k(s)$, the fraction field of $k[s]$), the ideal $(x^{pq},y^{pq}, g)C$ is primary to $P' = (x,y)C$, which is a maximal ideal of $C$.  So to show that $((x^q, y^q)^p + (g))C_{P'} / ((x^{pq},y^{pq}) +(g))C_{P'} = ((x^q, y^q)^p + (g) / (x^{pq}, y^{pq},g))_{P'}$ is nonzero over $A_P = C_{P'}$, it suffices to show that the $C$-module $((x^q, y^q)^p + (g))C / (x^{pq},y^{pq}, g)C \neq 0$.  For this, it is enough to show that $x^q y^{(p-1)q} \notin \ic := (x^{pq},y^{pq}, xy(x-y))C$, since $xy(x-y)$ is a factor of $g$.  Suppose that $x^q y^{(p-1)q} \in \ic$.  Then there exist polynomials $a,b,c \in C$ such that \[
x^q y^{(p-1)q} = a x^{pq} + b y^{pq} + cxy(x-y).
\]
From degree considerations (taking the homogeneous degree $pq$-part of the above equation), we may assume that $a,b \in L$.  Then making the substitution ($x=1$, $y=0$) in the displayed equation yields $a=0$, whereas the substitution ($x=0$, $y=1$) yields $b=0$.   So $x^q y^{(p-1)q} = cxy(x-y)$.  But then the substitution $x=y=1$ leads to the conclusion that $1=0$, a manifest contradiction.  Hence $x^q y^{(p-1)q} \notin \ic$, so that  $\p \in \Ass_R(I^{[q]}/ J^{[q]})$, as required.

We also know from Construction~\ref{cons} that $\m = (J^{[q]} : z)$, so that since $z\in I^{[q]}$, we have $\m \in \Ass_R(I^{[q]}/J^{[q]})$ as well.  So we are in the ``critical situation'' described at the beginning of this note.

Moreover, \[
\len_R(H^0_\m(I^{[q]}/J^{[q]})) \leq \len_A(H^0_\m(A/\ie)) = \len_A(A/\m) = 1,
\]
a constant, which shows that $\rjj{R_\m}{J_\m}{I_\m} = 0$, since $\dim R/J = 1>0$.  Hence, the expected converse to Theorem~\ref{thm:zerotc} holds for this specific example.
\end{example}

\section*{Acknowlegements}
We wish to thank Ezra Miller and Kirsten Schmitz for discussions regarding 
this note. We used Macaulay 2 \cite{M2hyper} for some of the computations.

\providecommand{\bysame}{\leavevmode\hbox to3em{\hrulefill}\thinspace}
\providecommand{\MR}{\relax\ifhmode\unskip\space\fi MR }
\providecommand{\MRhref}[2]{%
  \href{http://www.ams.org/mathscinet-getitem?mr=#1}{#2}
}
\providecommand{\href}[2]{#2}

\end{document}